\documentclass[leqno,12pt]{amsart}
\usepackage{amssymb} 
\theoremstyle{plain} 
\newtheorem{theorem}{\indent\sc Theorem}[section] 

\newtheorem{proposition}[theorem]{\indent\sc Proposition}

\theoremstyle{definition} 

\begin{document}

\title[Coupled Painlev\'e VI system]{Coupled Painlev\'e VI systems in dimension four with affine Weyl group symmetry of type $E_6^{(2)}$ \\}
\author{Yusuke Sasano }

\renewcommand{\thefootnote}{\fnsymbol{footnote}}
\footnote[0]{2000\textit{ Mathematics Subjet Classification}.
34M55; 34M45; 58F05; 32S65.}

\maketitle

\begin{abstract}
We find a four-parameter family of coupled Painlev\'e VI systems in dimension four with affine Weyl group symmetry of type $E_6^{(2)}$. This is the first example which gave higher order Painlev\'e type systems of type $E_{6}^{(2)}$. We study its symmetry and holomorphy conditions.
\end{abstract}

\section{Introduction}
In \cite{Sasa1,Sasa2,Sasa3,Sasa4,Sasa5,Sasa6,Sasa7}, we presented some types of coupled Painlev\'e systems with various affine Weyl group symmetries by connecting the invariant divisors $p_i,q_i-q_{i+1},p_{i+1}$ for the canonical variables $(q_i,p_i) \ (i=1,2,\ldots,n)$. These systems are polynomial Hamiltonian systems with coupled Painlev\'e Hamiltonians.

In this paper, we find a 4-parameter family of coupled Painlev\'e VI systems in dimension four with affine Weyl group symmetry of type $E_6^{(2)}$ given by
\begin{align}\label{1}
\begin{split}
\frac{dx}{dt}&=\frac{\partial H}{\partial y}, \quad \frac{dy}{dt}=-\frac{\partial H}{\partial x}, \quad \frac{dz}{dt}=\frac{\partial H}{\partial w}, \quad \frac{dw}{dt}=-\frac{\partial H}{\partial z}
\end{split}
\end{align}
with the polynomial Hamiltonian
\begin{align}\label{2}
\begin{split}
t(1-t)H =&x^2y^3+((1-2t)x-2\alpha_1-\alpha_2-\alpha_3)xy^2\\
&+\{(t-1)tx^2+((4\alpha_1+4\alpha_2+3\alpha_3+\alpha_4)t-(2\alpha_1+2\alpha_2+2\alpha_3+\alpha_4))x\\
&+\alpha_1(\alpha_1+\alpha_2+\alpha_3)\}y-(1-t)t \alpha_0 x\\
&+\frac{1}{4}[-z^2w^4+2 \alpha_3 zw^3+((1+t)z^2+2(2\alpha_1+2\alpha_2+\alpha_3)z- \alpha_3^2)w^2\\
&-2\{((-2\alpha_1-2\alpha_2-\alpha_3-\alpha_4)t+(2\alpha_1+2\alpha_2+2\alpha_3+\alpha_4))z\\
&+\alpha_3(2\alpha_1+2\alpha_2+\alpha_3)\}w-t(z+4\alpha_1+4\alpha_2+2\alpha_3)z]\\
&+(txz+(1-t)xzw-xzw^2-xyz+xyzw-\alpha_1(w-1)z+\alpha_3 xw)y.
\end{split}
\end{align}
Here $x,y,z$ and $w$ denote unknown complex variables, and $\alpha_0,\alpha_1, \dots ,\alpha_4$ are complex parameters satisfying the relation:
\begin{equation}\label{3}
\alpha_0+2\alpha_1+3\alpha_2+2\alpha_3+\alpha_4=1.
\end{equation}
In section 2, each principal part of this Hamiltonian can be transformed into canonical Painlev\'e VI Hamiltonian \eqref{HPVI} by birational and symplectic transformations.

This is the first example which gave higher order Painlev\'e type systems of type $E_{6}^{(2)}$.

We remark that for this system we tried to seek its first integrals of polynomial type with respect to $x,y,z,w$. However, we can not find. Of course, the Hamiltonian $H$ is not the first integral.

It is known that the Painlev\'e VI system admits the affine Weyl group symmetry of type $F_4^{(1)}$ (see \cite{Oka1}) as the group of its B{\"a}cklund transformations in addition to the diagram automorphisms of type $D_4^{(1)}$. The diagram automorphisms change the time variable $t$. However, in section 3, the system \eqref{1} admits the affine Weyl group symmetry of type $E_6^{(2)}$ as the group of its B{\"a}cklund transformations, whose generators $s_0,s_1,\ldots,s_4$ are determined by the invariant divisors \eqref{invariant}. Of course, these transformations do not change the time variable $t$.

\section{Principal parts of the Hamiltonian}
In this section, we study two Hamiltonians $K_1$ and $K_2$ in the Hamiltonian $H$.

At first, we study the Hamiltonian system
\begin{align}
\begin{split}
\frac{dx}{dt}&=\frac{\partial K_1}{\partial y}, \quad \frac{dy}{dt}=-\frac{\partial K_1}{\partial x}
\end{split}
\end{align}
with the polynomial Hamiltonian
\begin{align}\label{4}
\begin{split}
t(1-t)K_1 =&x^2y^3+((1-2t)x-2\alpha_1-\alpha_2-\alpha_3)xy^2\\
&+\{(t-1)tx^2+((4\alpha_1+4\alpha_2+3\alpha_3+\alpha_4)t-(2\alpha_1+2\alpha_2+2\alpha_3+\alpha_4))x\\
&+\alpha_1(\alpha_1+\alpha_2+\alpha_3)\}y-(1-t)t \alpha_0 x,
\end{split}
\end{align}
where setting $z=w=0$ in the Hamiltonian $H$, we obtain $K_1$.

We transform the Hamiltonian \eqref{4} into the Painlev\'e VI Hamiltonian:
\begin{align}\label{HPVI}
\begin{split}
&H_{VI}(x,y,t;\beta_0,\beta_1,\beta_2,\beta_3,\beta_4)\\
&=\frac{1}{t(t-1)}[y^2(x-t)(x-1)x-\{(\beta_0-1)(x-1)x+\beta_3(x-t)x\\
&+\beta_4(x-t)(x-1)\}y+\beta_2(\beta_1+\beta_2)x]  \quad (\beta_0+\beta_1+2\beta_2+\beta_3+\beta_4=1). 
\end{split}
\end{align}

{\bf Step 1:} We make the change of variables:
\begin{equation}
x_1=x, \quad y_1=y-t.
\end{equation}

{\bf Step 2:} We make the change of variables:
\begin{equation}
x_2=-y_1, \quad y_2=x_1.
\end{equation}
Then, we can obtain the Painlev\'e VI Hamiltonian:
\begin{align}
\begin{split}
&H_{VI}(x_2,y_2,t;\alpha_0,\alpha_2+\alpha_3,\alpha_1,\alpha_2+\alpha_3+\alpha_4,\alpha_2).
\end{split}
\end{align}
Of course, the parameters $\alpha_i$ and $\beta_j$ satisfy the relations:
\begin{equation}
\beta_0+\beta_1+2\beta_2+\beta_3+\beta_4=\alpha_0+2\alpha_1+3\alpha_2+2\alpha_3+\alpha_4=1.
\end{equation}
We remark that all transformations are symplectic.

Next, we study the Hamiltonian system
\begin{align}
\begin{split}
\frac{dx}{dt}&=\frac{\partial K_2}{\partial y}, \quad \frac{dy}{dt}=-\frac{\partial K_2}{\partial x}
\end{split}
\end{align}
with the polynomial Hamiltonian
\begin{align}\label{5}
\begin{split}
t(1-t)K_2 =&\frac{1}{4}[-z^2w^4+2 \alpha_3 zw^3+((1+t)z^2+2(2\alpha_1+2\alpha_2+\alpha_3)z- \alpha_3^2)w^2\\
&-2\{((-2\alpha_1-2\alpha_2-\alpha_3-\alpha_4)t+(2\alpha_1+2\alpha_2+2\alpha_3+\alpha_4))z\\
&+\alpha_3(2\alpha_1+2\alpha_2+\alpha_3)\}w-t(z+4\alpha_1+4\alpha_2+2\alpha_3)z],
\end{split}
\end{align}
where setting $x=y=0$ in the Hamiltonian $H$, we obtain $K_2$.

Let us transform the Hamiltonian \eqref{5} into the Painlev\'e VI Hamiltonian.

{\bf Step 1:} We make the change of variables:
\begin{equation}
t=T_1^2.
\end{equation}
We note that
\begin{equation}
dK_2 \wedge dt=2T_1d{\tilde K}_2 \wedge dT_1.
\end{equation}

{\bf Step 2:} We make the change of variables:
\begin{equation}
z_1=2z, \quad w_1=\frac{1}{2}w+\frac{1}{2}.
\end{equation}
By this transformation, in the coordinate system $(Z_1,W_1)=(1/z_1,w_1)$  two of four accessible singular points are transformed into $W_1=0$ and $W_1=1$.

{\bf Step 3:} We make the change of variables:
\begin{equation}
z_2=-(z_1w_1-\alpha_3)w_1, \quad w_2=\frac{1}{w_1}.
\end{equation}

{\bf Step 4:} We make the change of variables:
\begin{equation}
z_3=\frac{T_1-1}{T_1+1}z_2, \quad w_3=\frac{T_1+1}{T_1-1}w_2+\frac{2}{1-T_1}, \quad T_1=1-2T_2+2\sqrt{T_2(T_2-1)} .
\end{equation}
By this transformation, in the coordinate system $(Z_2,W_2)=(1/z_4,w_4)$  the others are transformed into $W_2=0$ and $W_2=\frac{1}{T_2}$. We remark that it is not $W_2 = \infty$ but  $W_2=0$ because we consider in the coordinate system $(z_2,w_2)$.

{\bf Step 5:} We make the change of variables:
\begin{equation}
z_4=-(z_3w_3-\alpha_3)w_3, \quad w_4=\frac{1}{w_3}.
\end{equation}

{\bf Step 6:} We make the change of variables:
\begin{equation}
z_5=w_4, \quad w_5=-z_4.
\end{equation}
Then, we can obtain the Painlev\'e VI Hamiltonian:
\begin{align}
\begin{split}
&\frac{1}{2} H_{VI}(z_5,w_5,T_2;\alpha_0+\alpha_2-1,\alpha_0+\alpha_2,\alpha_3,\alpha_4,1-\alpha_0+2\alpha_1+\alpha_2).
\end{split}
\end{align}
We remark that all transformations are symplectic.

\section{Symmetry and holomorphy conditions}
In this section, we study the symmetry and holomorphy conditions of the system \eqref{1}. These properties are new.

\begin{theorem}\label{th:1}
The system \eqref{1} admits the affine Weyl group symmetry of type $E_6^{(2)}$ as the group of its B{\"a}cklund transformations, whose generators $s_0,s_1,\ldots,s_4$ defined as follows$:$ with {\it the notation} $(*):=(x,y,z,w,t;\alpha_0,\alpha_1,\ldots,\alpha_4)$\rm{: \rm}

\begin{figure}
\unitlength 0.1in
\begin{picture}( 54.3300,  8.4600)( 14.4100,-16.0600)
%
\special{pn 20}%
\special{ar 1858 1188 418 418  0.0000000 6.2831853}%
%
\special{pn 20}%
\special{ar 2968 1188 418 418  0.0000000 6.2831853}%
%
\special{pn 20}%
\special{ar 4128 1178 418 418  0.0000000 6.2831853}%
%
\special{pn 20}%
\special{ar 5288 1178 418 418  0.0000000 6.2831853}%
%
\special{pn 20}%
\special{ar 6458 1188 418 418  0.0000000 6.2831853}%
\put(17.2000,-12.6000){\makebox(0,0)[lb]{$y$}}%
\put(28.1000,-12.5200){\makebox(0,0)[lb]{$x$}}%
\put(37.4100,-12.7200){\makebox(0,0)[lb]{$y+w^2-t$}}%
\put(51.1700,-12.5600){\makebox(0,0)[lb]{$z$}}%
\put(62.1900,-12.6600){\makebox(0,0)[lb]{$w-1$}}%
%
\special{pn 20}%
\special{pa 2290 1190}%
\special{pa 2550 1190}%
\special{fp}%
%
\special{pn 20}%
\special{pa 3410 1190}%
\special{pa 3690 1190}%
\special{fp}%
%
\special{pn 20}%
\special{pa 5720 1190}%
\special{pa 6020 1190}%
\special{fp}%
%
\special{pn 20}%
\special{pa 4520 1020}%
\special{pa 4910 1020}%
\special{fp}%
\special{sh 1}%
\special{pa 4910 1020}%
\special{pa 4844 1000}%
\special{pa 4858 1020}%
\special{pa 4844 1040}%
\special{pa 4910 1020}%
\special{fp}%
%
\special{pn 20}%
\special{pa 4530 1350}%
\special{pa 4900 1350}%
\special{fp}%
\special{sh 1}%
\special{pa 4900 1350}%
\special{pa 4834 1330}%
\special{pa 4848 1350}%
\special{pa 4834 1370}%
\special{pa 4900 1350}%
\special{fp}%
\end{picture}%
\label{fig:E6Painlevefig1}
\caption{The figure denotes the Dynkin diagram of type $E_6^{(2)}$. The symbol in each circle denotes the invariant divisors of the system \eqref{1} of type $E_6^{(2)}$.}
\end{figure}
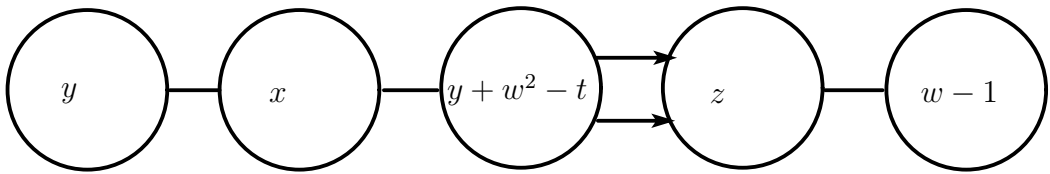
\begin{align}
\begin{split}
s_0:(*) \rightarrow &\left(x+\frac{\alpha_0}{y},y,z,w,t;-\alpha_0,\alpha_1+\alpha_0,\alpha_2,\alpha_3,\alpha_4 \right),\\
s_1:(*) \rightarrow &\left(x,y-\frac{\alpha_1}{x},z,w,t;\alpha_0+\alpha_1,-\alpha_1,\alpha_2+\alpha_1,\alpha_3,\alpha_4 \right),\\
s_2:(*) \rightarrow &\left(x+\frac{\alpha_2}{y+w^2-t},y,z+\frac{2\alpha_2 w}{y+w^2-t},w,t;\alpha_0,\alpha_1+\alpha_2,-\alpha_2,\alpha_3+2\alpha_2,\alpha_4 \right),\\
s_3:(*) \rightarrow &\left(x,y,z,w-\frac{\alpha_3}{z},t;\alpha_0,\alpha_1,\alpha_2+\alpha_3,-\alpha_3,\alpha_4+\alpha_3 \right),\\
s_4:(*) \rightarrow &\left(x,y,z+\frac{\alpha_4}{w-1},w,t;\alpha_0,\alpha_1,\alpha_2,\alpha_3+\alpha_4,-\alpha_4 \right).
\end{split}
\end{align}
\end{theorem}
We note that the B{\"a}cklund transformations of this system satisfy
\begin{equation}
s_i(g)=g+\frac{\alpha_i}{f_i}\{f_i,g\}+\frac{1}{2!} \left(\frac{\alpha_i}{f_i} \right)^2 \{f_i,\{f_i,g\} \}+\cdots \quad (g \in {\Bbb C}(t)[x,y,z,w]),
\end{equation}
where poisson bracket $\{,\}$ satisfies the relations:
$$
\{y,x\}=\{w,z\}=1, \quad the \ others \ are \ 0.
$$
Since these B{\"a}cklund transformations have Lie theoretic origin, similarity reduction of a Drinfeld-Sokolov hierarchy admits such a B{\"a}cklund symmetry.

\begin{proposition}
This system has the following invariant divisors\rm{:\rm}
\begin{center}\label{invariant}
\begin{tabular}{|c|c|c|} \hline
parameter's relation & $f_i$ \\ \hline
$\alpha_0=0$ & $f_0:=y$  \\ \hline
$\alpha_1=0$ & $f_1:=x$  \\ \hline
$\alpha_2=0$ & $f_2:=y+w^2-t$  \\ \hline
$\alpha_3=0$ & $f_3:=z$  \\ \hline
$\alpha_4=0$ & $f_4:=w-1$  \\ \hline
\end{tabular}
\end{center}
\end{proposition}
We note that when $\alpha_0=0$, we see that the system \eqref{1} admits a particular solution $y=0$, and when $\alpha_2=0$, after we make the birational and symplectic transformations:
\begin{equation}
x_2=x, \ y_2=y+w^2-t, \ z_2=z-2xw, \ w_2=w
\end{equation}
we see that the system \eqref{1} admits a particular solution $y_2=0$.

\begin{theorem}\label{pro:2}
Let us consider a polynomial Hamiltonian system with Hamiltonian $K \in {\Bbb C}(t)[x,y,z,w]$. We assume that

$(A1)$ $deg(K)=6$ with respect to $x,y,z,w$.

$(A2)$ This system becomes again a polynomial Hamiltonian system in each coordinate system $r_i \ (i=0,1,\ldots,4)${\rm : \rm}
\begin{align}
\begin{split}
r_0:&x_0=\frac{1}{x}, \ y_0=-(yx+\alpha_0)x, \ z_0=z, \ w_0=w, \\
r_1:&x_1=-(xy-\alpha_1)y, \ y_1=\frac{1}{y}, \ z_1=z, \ w_1=w,\\
r_2:&x_2=\frac{1}{x}, \ y_2=-\left((y+w^2-t)x+\alpha_2 \right)x, \ z_2=z-2xw, \ w_2=w, \\
r_3:&x_3=x, \ y_3=y, \ z_3=-(zw-\alpha_3)w, \ w_3=\frac{1}{w}, \\
r_4:&x_4=x, \ y_4=y, \ z_4=\frac{1}{z}, \ w_4=-((w-1)z+\alpha_4)z.
\end{split}
\end{align}
Then such a system coincides with the system \eqref{1} with the polynomial Hamiltonian  \eqref{2}.
\end{theorem}
By this theorem, we can also recover the parameter's relation \eqref{3}.

We note that the condition $(A2)$ should be read that
\begin{align*}
&r_j(K) \quad (j=0,1,3,4), \quad r_2(K+x)
\end{align*}
are polynomials with respect to $x,y,z,w$.

\end{document}